\documentclass[sn-mathphys-num]{sn-jnl}

\usepackage{graphicx}%
\usepackage{multirow}%
\usepackage{amsmath,amssymb,amsfonts}%
\usepackage{amsthm}%
\usepackage{mathrsfs}%
\usepackage[title]{appendix}%
\usepackage{xcolor}%
\usepackage{textcomp}%
\usepackage{manyfoot}%
\usepackage{booktabs}%
\usepackage{algorithm}%
\usepackage{algorithmicx}%
\usepackage{algpseudocode}%
\usepackage{listings}%

\newtheorem{theorem}{Theorem}
\newtheorem{proposition}[theorem]{Proposition}%
\newtheorem{lemma}[theorem]{Lemma}

\newtheorem{remark}{Remark}%

\def\R{\mathbb{R}}

\begin{document}

\title[Article Title]{Local controllability of free boundary three-dimensional semilinear radial parabolic equations}


\author*[1]{\fnm{Juan} \sur{L\'imaco}}\email{jlimaco@id.uff.br}

\author[1,2]{\fnm{Luis P.} \sur{Yapu}}


\affil*[1]{\orgdiv{Instituto de Matem\'atica e Estat\'\i stica}, \orgname{Universidade Federal Fluminense}, \orgaddress{\street{Rua Prof. Marcos Waldemar de Freitas Reis S/N (Campus do Gragoat\'a)}, \city{Niteroi}, \postcode{CEP 24210-201}, \state{Rio de Janeiro}, \country{Brazil}}}
\affil*[2]{\orgdiv{Chair of Dynamics, Control, Machine Learning and Numerics}, \orgname{Friedrich-Alexander Universität Erlangen-Nürnberg (FAU)},
\orgaddress{\city{Erlangen}, \country{Germany}}}




\abstract{We prove that a free boundary semilinear heat equation with Stefan boundary condition and radially symmetric data is locally null controllable. 
The strategy involves reducing the problem to the corresponding one-dimensional formulation and adapting a Carleman inequality in that setting.
The local null controllability of the free-boundary problem is then established via the Schauder fixed-point theorem.
To the best of our knowledge, this is the first controllability result for this problem with Stefan boundary condition in more than one spatial dimension.}




\maketitle

\section{Introduction} \label{sec:Intro}

The free-boundary problem with the Stefan boundary condition is a classical problem about the evolution of the temperature of a system composed by phases of a fluid. More specifically, in the heat evolution model that we consider, an evolving phase composed of water is surrounded by ice at zero temperature. Since the ice temperature does not evolve, this is called a one-phase Stefan problem.
The equation models the evolution of the temperature $y(x,t)$ at some position and time. The water phase is in the time-dependent domain $\Omega(t) \subset \R^3$, with boundary $\partial \Omega(t)$ that separates the water and the ice.

\subsection{Main contributions}
In this paper, we are interested in the case of a three-dimensional one-phase Stefan problem with radial symmetry.
Let $B_R$ denote the Euclidean three-dimensional ball of radius $R>0$ centered at the origin.  
Consider the domain $Q_R = \{(x,t) : x \in B_{R(t)}, t \in (0,T) \} \subset \R^3 \times \R$, with lateral boundary $\Sigma_R = \{(x,t) : x \in \partial B_{R(t)}, t \in (0,T) \}$, where $R(t) \in C^1([0,T])$ describes the radius that changes with time and verifies, for some positive constants $R_*$ and $E$, 
\begin{equation}
\label{eq:bound_on_R}
0 < R_* \leq R(t) \leq E.
\end{equation}

Let $f : \R \to \R$ be a globally Lipschitz continuous function such that $f(0)=0$.
Fix the control region $\omega := B_b$, where $0 < b < R_*$ and consider the controlled free-boundary problem
\begin{equation}
\label{eq:PDE-intro}
\begin{cases}
 y_t - \Delta y + f(y) = 1_{\omega} u, & \text{in } Q_R, \\
y(x,t) = 0, & \text{on } \Sigma_R, \\
 y(\cdot,0) = y_0, & \text{in } B_{R(0)},
\end{cases}
\end{equation}
where $u(x,t)$ indicates the control and the Stefan condition describing the evolution of the boundary $\partial B_{R(t)}$ is given by
\begin{equation}
\label{eq:stefan_condition-intro}
R'(t) = -\nabla y(x,t) \cdot \hat n, \quad x \in \partial B_{R(t)},
\end{equation}
where $\hat n$ denotes the outward-pointing unit normal vector at $(x,t)$.


System \eqref{eq:PDE-intro} is said \emph{locally null controllable} at time $T>0$ if there exists a constant $\varepsilon>0$ such that, if $\|y_0\|_{H_0^1(B_{R(0)})} \leq \varepsilon$, there is a solution $(y,u,R)$ of \eqref{eq:PDE-intro}, \eqref{eq:stefan_condition-intro}, such that 
\begin{align*}
R \in C^1([0,T]), \qquad & 0 < R_* \leq R(t) \leq E, \\
u \in L^2(\omega \times (0,T)), \qquad & y \in C^0([0,T],H_0^1(B_{R(t)}),
\end{align*}
and $y(x,T) = 0$, for $x \in B_{R(T)}$.

The main result is the following local null controllability theorem in the three-dimensional radially symmetric case: 

\begin{theorem}
\label{thm:main_theorem}
Assume that there are constants $0 < b < R_* < R_0 < E$.
Let $T>0$ and $y_0 \in H_0^1(B_{R_0})$ radially symmetric. 
Then, putting $R(0)=R_0$, the three-dimensional free-boundary problem \eqref{eq:PDE-intro}-\eqref{eq:stefan_condition-intro} is locally null controllable.
\end{theorem}

The result is a consequence of the controllability of the associated linearized system and the use of Schauder's fixed-point theorem to the nonlinear system.
This is the first result of controllability for the free-boundary problem \eqref{eq:PDE-intro} in more than one spatial dimension using the Stefan condition \eqref{eq:stefan_condition-intro}, rather than a modified condition such as the one considered in \cite{Demarque_EFC-18} in the two-dimensional case and star-shaped domain.

Although the Stefan problem has been extensively studied \cite{Aronson-83, Friedman-82, Meirmanov-92}, results concerning its controllability are recent. Most known controllability results pertain to the one-dimensional case.

\subsection{Related work}

The first local null controllability result for the heat equation with free boundary in one dimension was given by Fernández-Cara, Límaco and Menezes in \cite{EFC_Lim_Men-16}. Soon after, this was generalized to include a semilinearity by Fernández-Cara and Sousa in \cite{EFC_Ivaldo-17}. Both papers treat the intrinsic nonlinearity of the problem using a fixed-point strategy and are limited to the one-dimensional case because the fixed point map uses an explicit integration of the one-dimensional Stefan condition. Afterwards, Fernández-Cara, Hernández and Límaco in \cite{EFC_Lim_Freddy-18} proved the same results using an inverse function theorem, which could in principle be suitable to treat the higher dimensional case, but no further results where published in that direction. The only related work we are aware of in two dimensions is by Demarque and Fernández-Cara \cite{Demarque_EFC-18}, where the condition describing the evolution of the free boundary was modified introducing a diffusion term.

In \cite{EFC_Lim_Freddy-18} the authors used a compactness-uniqueness method in order to treat de Stefan condition as a constraint. Geshkovski and Zuazua used the same method to treat the Burgers equation with free boundary in \cite{Gesh_Zua-21}. Moreover, Araujo, Fernández-Cara, Límaco and Souza in \cite{Arau_EFC_Sou-22} used also compactness-uniqueness arguments to prove the local null controllability of the two-phase Stefan problem and also to control the boundary curve between the two phases.

More recently, Barcena-Petisco, Fernández-Cara and Souza \cite{Bar_EFC_Sou-23} proved the controllability to trajectories of the one-phase Stefan problem. Also Costa,
Límaco, Lopes and Prouvée \cite{Vitor_JAL-23} included in the model nonlocal nonlinearities. This required a more careful study and Hölder regularity of the initial state. Moreover,  Wang, Lan, Lei in \cite{WLL-22} studied the quasilinear case, which also needs Hölder regularity of the initial state. On the other hand, Wang, Lei and Wu \cite{WLW-22} studied the existence of insensitizing controls for the semilinear Heat equation.

In a different geometry, Geshkovski and Maity in \cite{Gesh_Maity-23} studied the Stefan problem in a periodic box. Their methods are different, using the periodicity to find a Fourier decomposition in one dimension and the study of the resultant sequence of one-dimensional systems and proving an observability result which is uniform for the sequence.

\subsection{Overview}

This paper is organized as follows. In Section \ref{sec:reduction} we present the reduction of the radial problem to the one-dimensional case and show the null controllability of the linearized system using an adapted Carleman inequality. In Section \ref{sec:proof_main} we prove
Theorem \ref{thm:main_theorem} using Schauder's fixed point theorem. Finally, in section \ref{sec:additional_remarks} we observe that boundary controllability holds directly when both the domain and the control regions are annular (excluding the origin). We also present several problems that are already known to be locally null controllable in the one-dimensional case.

\section{Controllability of the associated one-dimensional  linearized system}
\label{sec:reduction}

In this section, using the radial symmetry we reduce the equation to one dimension and after a change of variables, we prove the controllability of the system without the Stefan condition, supposing that the boundary follows a path determined by a function $R(t) \in C^1([0,T])$ such that \eqref{eq:bound_on_R} holds.



Using the radial symmetry, the Laplacian term can be written in polar coordinates. We get an equation applied to the radial function  $z(r,t) = y(x,t)$, with $r=|x|$. In this case, the 3D Laplacian is given by $\Delta z =  \frac{1}{r^2} (r^2 z_r)_r = z_{rr} + \frac{2}{r} z_r$, with radial control $v(r,t)=u(x,t)$, and the system \eqref{eq:PDE-intro} becomes
\begin{equation}
\label{eq:1D_PDE_semilin}
\begin{cases}
z_t - z_{rr} - \frac{2}{r} z_r + f(z) = 1_{\omega} v, & \text{in } \hat Q_R, \\
z_r (0,t) = 0, \quad z(R(t),t) = 0, & \text{in } (0,T),\\
z(\cdot,0) = z_0, & \text{in } (0,R(0)),
\end{cases}
\end{equation}
where the domain $\hat Q_R = \{(r,t) : r \in [0,R(t)], \ t \in (0,T) \}$ is one-dimensional in space. Let us observe that $y(\cdot,t) \in L^2(B_{R(t)})$ for a.e. $t \in (0,T)$ is equivalent to $z(\cdot,t) \in L^2([0,R(t)], r^2 dr)$, where in the first case the $L^2$ space is in three dimensions and in the second case it is a weighted $L^2$ space in one dimension, i.e., for a.e. $t \in (0,T)$,  we have 
$\int_0^{R(t)} |z(r,t)|^2 r^2 dr < \infty$. On the other hand, 
the Stefan condition \eqref{eq:stefan_condition-intro} becomes
\begin{equation}\label{eq:stefan_condition_1}
    R'(t) = -z_r(R(t),t), \qquad t \in (0,T).
\end{equation}

The linearized system about zero (stationary solution of \eqref{eq:1D_PDE_semilin}) is given by 
\begin{equation}
\label{eq:1D-PDE-linear}
\begin{cases}
 z_t - z_{rr} - \frac{2}{r} z_r + a(r,t) z = 1_{\omega} v(r,t), & \text{in}\quad \hat Q_R, \\
z_r (0,t) = 0, \qquad z(R(t),t) = 0, & \text{in} \quad (0,T),\\
z(\cdot,0) = z_0, & \text{in}\quad (0,R(0)),
\end{cases}
\end{equation}
where $a(r,t)$ is a bounded function.

In order to remove the singular term at $r=0$, we define the functions $\tilde z(r,t) = r \, z(r,t)$ and $\tilde v = r \, v(r,t)$. Then, we have that $\tilde z(0,t) = 0$ and \eqref{eq:1D-PDE-linear} becomes
\begin{equation}
\label{eq:1D-PDE-change}
\begin{cases}
\tilde z_t - \tilde z_{rr} + a(r,t) \tilde z = 1_{\omega} \tilde v(r,t), & \text{in}\quad \hat Q_R, \\
\tilde z (0,t) = 0, \quad \tilde z(R(t),t) = 0, & \text{in}\quad (0,T), \\
\tilde z(\cdot,0) = \tilde z_0, & \text{in}\quad (0,R(0)).
\end{cases}
\end{equation}

Note that after the change of variable, $z(\cdot,t) \in L^2([0,R(t)], r^2 dr)$ is equivalent to $\tilde z(\cdot,t) \in L^2([0,R(t)],dr)$.

The Stefan condition \eqref{eq:stefan_condition_1} is transformed into
\begin{equation}
\label{eq:stefan_condition_2}
R'(t) = \frac{\tilde z_r(R(t),t)}{R(t)} - \frac{\tilde z(R(t),t)}{R(t)^2} = \frac{\tilde z_r(R(t),t)}{R(t)},
\end{equation}
where in the last equality we used the boundary condition. We observe, from the condition $0 < R_* < R(t) < E$, that the denominator is not zero.

\begin{remark}[On the $N$-dimensional case]
The expression of the Laplacian used above is valid in tree spatial dimensions. If the system is $N$-dimensional, then the radial Laplacian is given by $\Delta z =  z_{rr} + \frac{N-1}{r} z_r$ and after the change of variable $\tilde z(r,t) = r^{\frac{N-1}{2}} z(r,t)$, the equation analogous to \eqref{eq:1D-PDE-change} becomes
$$\tilde z_t - \tilde z_{rr} + \frac{(N-1)(N-3)}{4r^2} \tilde z + a(r,t) \tilde z = 1_\omega \tilde v(r,t),$$
which still has a singularity at the origin if $N \neq 1$ or $N \neq 3$. That kind of system has been studied in \cite{Vancos-Zuazua-08} but not in the free-boundary case.  
\end{remark}

The adjoint system of \eqref{eq:1D-PDE-change} is
\begin{equation}
\label{eq:1D-PDE-adjoint}
\begin{cases}
-\phi_t - \phi_{rr}  + a(r,t) \phi = F(r,t), & \text{in}\quad \hat Q_R, \\
\phi (0,t) = 0, \quad \phi(R(t),t) = 0, & \text{on}\quad (0,T),\\
\phi(\cdot,T) = \phi_T. & \text{in}\quad (0,R(T)),
\end{cases}
\end{equation}
where $F \in L^2(\hat Q_R)$ and $\phi_T \in L^2([0,R(T)])$.

The adjoint system \eqref{eq:1D-PDE-adjoint} has the same form as in the one dimensional case studied in \cite{EFC_Ivaldo-17}. Analogously, the following observability property is a consequence of the adapted Carleman inequality of Section \ref{subsec:carleman}: 

\begin{proposition}\label{prop:observability}
There exists a constant $C > 0$, such that for any $\phi_T \in L^2([0,R(T)])$, the solution of \eqref{eq:1D-PDE-adjoint} with $F = 0$ satisfies
$$\int_0^{R(0)} |\phi(r,0)|^2 dr \leq C \iint_{[0,b) \times (0,T)}  |\phi|^2 dr dt.$$
\end{proposition}

This observability property implies an approximate controllability result of \eqref{eq:1D-PDE-change} by a standard argument of Fabre-Puel-Zuazua \cite{Fabre_Puel_Zuazua-95} using minimizers of a functional $J_\epsilon(\cdot;g,R)$, see formula \eqref{eq:J_eps_functional} and more details in Section 2.3 of \cite{EFC_Lim_Men-16}. By uniformity with respecto to $\epsilon$, the approximate controllability implies the null controllability result of the linearized system.

\begin{theorem}
For any $\tilde z_0 \in H_0^1([0,R(0)])$ there exist a pair $(\tilde v, \tilde z)$, with $\tilde v \in L^2([0,b) \times (0,T))$ and $\tilde z \in C^0([0,T],H_0^1([0,R(T)]))$ solving \eqref{eq:1D-PDE-change} and \eqref{eq:stefan_condition_2} such that
$$\tilde z(x,T)=0, \qquad x \in (0,R(T)).$$
Moreover, $\tilde v$ verifies
\begin{equation}
\label{eq:v_bounded}
\|\tilde v\|_{L^2([0,b)\times (0,T))} \leq C \| \tilde z_0 \|_{H_0^1([0,R(0)])},
\end{equation}
where $C$ only depends on $R_*$, $E$, $[0,b)$, $\|R'\|_\infty$, $\|a\|_\infty$ and $T$.
\end{theorem}

We observe that the one-dimensional Carleman estimate used in  \cite{EFC_Ivaldo-17, EFC_Lim_Men-16} needs a modification for our problem because the radially symmetric control domain $\omega=B_b$ after passing to the one-dimensional problem becomes the interval $[0,b)$ which is not an open interval $(a,b)$
as in \cite{EFC_Ivaldo-17, EFC_Lim_Men-16}. 
This is fixed by taking the symmetric interval $(-R(t),R(t))$ and using an even extension of the weights in the Carleman estimate.

\subsection{Carleman estimate}
\label{subsec:carleman}

Let us set the notation $\hat Q_{RR} = \{(r,t) : r \in [-R(t),R(t)], \ t \in (0,T) \}$. We need the following lemma which is a simple adaptation of Lemma 2.1 in \cite{EFC_Lim_Men-16}.

\begin{lemma}
For some $b_0 < b$, there exists an even function $\eta_0 \in C^1(\hat Q_{RR})$, with $\eta_{0,rr} \in C^0(\hat Q_{RR})$, such that
\begin{equation*}
\begin{cases}
\alpha_0(r,t) = 0,\quad& (r,t) \in \{R(t),-R(t)\} \times (0,T) \\
|\alpha_{0,r}| > 0, \quad& (r,t) \in ((-R(t),-b_0) \cup (b_0,R(t)))\times (0,T) \\
\alpha_0(r,t) = 1 - \frac{r-b}{R(t)-b}, \quad& (r,t) \in ((-R(t),-b) \cup (b,R(t)))\times (0,T).
\end{cases} 
\end{equation*}
\end{lemma}

For the proof we can take, for instance, the even extension on $[-R(t),R(t)]$ of the following function, which is an adaptation from \cite{EFC_Lim_Men-16}, 
\begin{equation*}
\alpha_0(r,t) =
\begin{cases}
1 + p\left(\frac{2(b-r)}{2b}, \frac{2b}{2(R(t)-b)} \right), & \text{if} \quad 0 \leq r < b, \\
\frac{R(t)-r}{R(t)-b}, &\text{if}\quad b \leq r \leq R(t),
\end{cases}
\end{equation*}
where $p(w,z) = zw + (10-6z)w^3 + (8z-15)w^4 + (6-3z)w^5$. We verify easily that
$\alpha_{0,r}(0,t) = -\frac{1}{b} p_w(1,\frac{2b}{2(R(t)-b)})=0$.

For some parameter $\lambda >0$, we define 
$\alpha_1(r,t) = \alpha_0(r,t)+1$, 
$\sigma(r,t) = e^{2 \lambda \| \alpha_1 \|_\infty } - e^{\lambda\alpha_1(r,t)},$ 
and the weight functions, for $k \geq 2$,
$$\alpha (r,t) = \frac{\sigma(r,t)}{t^k(T-t)^k}, \qquad\qquad
\xi (r,t) = \frac{e^{\lambda \alpha_1(r,t) }}{t^k(T-t)^k}.$$

We use the notation
\begin{equation}
\label{eq:formula_I}
\begin{split}
I(\phi) := & \iint_{\hat Q_R} e^{-2s\alpha} \left( \frac{1}{s\xi} (|\phi_t|^2+|\phi_{rr}|^2) + \lambda^2 s\xi |\phi_r|^2 + \lambda^4 s^3 \xi^3 |\phi|^2 drdt \right)\\
& + \int_0^T e^{-2s\alpha(R(t),t)} \lambda s \xi(R(t),t) |\phi_r(R(t),t)|^2 dt. 
\end{split}    
\end{equation}

We obtained the following Carleman inequality:

\begin{proposition}
\label{prop:carleman}
There exist positive constants $C$, $s_1$ and $\lambda_1$ such that, for $s>s_1$, $\lambda > \lambda_1$, $F \in L^2(\hat Q_R)$, the solution of the adjoint system \eqref{eq:1D-PDE-adjoint} satisfies
\begin{equation}
\label{eq:carleman}
I(\phi) \leq C \left( \iint_{[0,b)\times(0,T)} \lambda^4 s^3 \xi^3 |\phi|^2 drdt + \iint_{\hat Q_R} e^{-2s\alpha} |F|^2 dr dt \right).
\end{equation}
\end{proposition}

The proof of this result can be found in Appendix \ref{appendix A} and this implies the observability inequality in Proposition \ref{prop:observability}.

\section{Local controllability of the nonlinear system}
\label{sec:proof_main}

In this section, we use Schauder's fixed-point theorem to show the local controllability of \eqref{eq:1D-PDE-change}. Substituting the original variables, we get the local null controllability of the original equation \eqref{eq:PDE-intro}.

The following compact embeddings from \cite{LSU-68} were already used in \cite{EFC_Ivaldo-17} in order to get compactness in the application of the Schauder fixed-point theorem.

\noindent\textbf{Hölder regularity of $\tilde z_r$:}

Given an integer $m\geq 0$, a constant $\alpha \in (0,1)$ and a non-empty open set $G$ in
$$\hat Q_0 = (0,E) \times (0,T),$$
we set
$$
\langle u \rangle_{r,G}^{(\alpha)} = \sup_{(r,t),(r',t) \in \bar G} \frac{|u(r,t)-u(r',t)|}{|r-r'|^\alpha},
\qquad
\langle u \rangle_{r,G}^{(m+\alpha)} = \sum_{2k+l=m} \langle D_t^k D_r^l u \rangle_{r,G}^{(\alpha)},
$$
$$
\langle u \rangle_{t,G}^{(\alpha/2)} = \sup_{(r,t),(r,t') \in \bar G} \frac{|u(r,t)-u(r,t')|}{|t-t'|^{\alpha/2}},
\qquad
\langle u \rangle_{t,G}^{\left(\frac{m+\alpha}{2} \right)} = \sum_{2k+l=m} \langle D_t^k D_r^l u \rangle_{r,G}^{(\alpha/2)},
$$
and
$$
|u|^{(m+\alpha)} := \sum_{2k+l \leq m} \|D_t^k D_r^l u\|_{L^\infty(G)} + \langle u \rangle_{r,G}^{(m+\alpha)} + \langle u \rangle_{t,G}^{\left(\frac{m+\alpha}{2} \right)}. 
$$

We denote by $K^{m,\alpha}$ the space of functions $u=u(r,t)$ such that $|u|_G^{(m+\alpha)}$ is finite. This is a Banach space with the norm $|u|_G^{(m+\alpha)}$. If $m+\alpha < m'+\alpha'$, we have the compact embedding $K^{m',\alpha'} \subset\subset K^{m,\alpha}$.

Let $V_R(t) = \tilde z_r(R(t),t)$, for any $t \in [0,T]$. From the estimates in \cite[Theorems 10.1-11.1]{LSU-68} we have that $\tilde z \in K^{1,\kappa}$ for any $\kappa \in [0,\frac{1}{2})$. Moreover, $V_R \in C^\kappa([0,T])$ and by the same estimates as in Section 2.3 of \cite{EFC_Ivaldo-17}, 
\begin{equation}
\label{eq:VR_bound}
\|V_R\|_{C^\kappa([0,T])} \leq C_1 \|\tilde z_0\|_{L^\infty([0,R(0)])}, 
\end{equation}
where the constant $C_1>0$ depends only on $\|a\|_{L^\infty(\hat Q_0)}$, $\|R'\|_\infty$, $T$, $R_*$, $[0,b)$ and $\|\tilde  z_0\|_{L^2([0,R(0)])}$.

\begin{proof}(of Theorem \ref{thm:main_theorem})

We follow closely the proof of the main theorem in \cite{EFC_Ivaldo-17} verifying the hypothesis of Schauder's fixed-point theorem.

First suppose that $f$ is of class $C^1$ and define the function $g:\R \to \R$ as $g(s)=\frac{f(s)}{s}$, if $s\neq 0 $ and $g(0)=f'(0)$. The case $f$ Lipschitz is obtained by approximation.

Let $(\bar z, \bar R) \in L^\infty(\hat Q_0) \times C^1([0,T])$ be a given pair of state and radius such that $R_* \leq \bar R(t) \leq E$. For any $\tilde z_0 \in H_0^1([0,R_0])$, $R_0 = \bar R(0)$, consider the linear approximate problem depending on $\epsilon>0$,
\begin{equation}
\label{eq:1D-PDE-aprox1}
\begin{cases}
\tilde z_{t,\epsilon} - \tilde z_{rr,\epsilon} + g(\bar z) \tilde z_\epsilon = 1_{\omega} \tilde v, & \text{in} \quad \hat Q_R, \\
\tilde z_\epsilon (0,t) = 0, \quad \tilde z_\epsilon(R(t),t) = 0, & \text{on} \quad (0,T), \\
\tilde z_\epsilon(\cdot,0) = \tilde z_0, & \text{in}\quad (0,R_0),
\end{cases}
\end{equation}
\begin{equation}
\label{eq:1D-PDE-aprox2}
\Vert \tilde z_\epsilon(\cdot, T) \Vert_{L^2([0,R(t)])} \leq \epsilon.
\end{equation}

Let us introduce the sets
\begin{equation}
\label{eq:espace_Z}
    \mathcal{Y} = \lbrace \bar z \in L^\infty(\hat Q_0) : \Vert \bar z \Vert_{L^\infty(\hat Q_0)} \leq K \rbrace,    
\end{equation}
and
$$\mathcal{Z} = \lbrace \bar R \in C^1([0,T]) : R_* \leq \bar R(t) \leq E, \bar R(0)=R_0, \Vert \bar R' \Vert_\infty \leq K_1 \rbrace, $$
where $K>0$ will be defined later and $K_1>0$ is given.

Consider the mapping
\begin{equation*}
\begin{split}
\Lambda_\epsilon : \mathcal{Y} \times \mathcal{Z} &\longrightarrow L^\infty (\hat Q_0) \times C^1([0,T]) \\
(\bar z, \bar R) &\mapsto (\tilde z_\epsilon, R_\epsilon),
\end{split}
\end{equation*}
where $\tilde z_\epsilon$ is a solution of the problem \eqref{eq:1D-PDE-aprox1}-\eqref{eq:1D-PDE-aprox2} for $\tilde v = \phi_\epsilon 1_\omega$, where $\phi_\epsilon$ is the solution of the adjoint system \eqref{eq:1D-PDE-adjoint} with final condition $\phi_T$  given by the unique minimizer $\phi_{\epsilon,T}$ of the functional 
\begin{equation}
\label{eq:J_eps_functional}
J_\epsilon(\phi_T;R(t)) = \frac{1}{2} \iint_{[0,b) \times (0,T)} |\phi|^2 drdt + \epsilon \|\phi_T\|_{L^2([0,R(T)])} + (\phi(\cdot,0),z_0)_{L^2([0,R_0])},
\end{equation}
for $z^0 \in L^2([0,R_0])$ and $\phi_T \in L^2([0,R(T)])$.

Moreover, integrating the Stefan condition \eqref{eq:stefan_condition_2} we get
\begin{equation}
\label{eq:stefan_integrado}
R_\epsilon(t) = R_0 - \int_0^t \frac{\tilde z_{\epsilon,r} (\bar R(s),s)}{\bar R(s)} ds.
\end{equation}

From the results of approximate controllability of the linear system \eqref{eq:1D-PDE-aprox1}-\eqref{eq:1D-PDE-aprox2}  
the map  $\Lambda_\epsilon$ is well defined. 

Now we show that $\Lambda_\epsilon$ maps $\mathcal{Y} \times \mathcal{Z}$ into itself. Indeed, we have that
\begin{equation}
\label{eq:z_eps-z0}
    \Vert \tilde z_\epsilon \Vert_{L^\infty(\hat Q_0)} \leq C_2 \Vert \tilde z_0 \Vert_{L^\infty([0,R_0])}, 
\end{equation}
where $C_2$ only depends on $R_*$, $E$, $[0,b)$, $K_1$ and $T$.

Moreover, since $R_* < R(t) < E$, using \eqref{eq:z_eps-z0} we get, from \eqref{eq:stefan_integrado},
$$
| R'_\epsilon(t) | \leq C_3 \Vert \tilde z_0 \Vert_{H_0^1([0,R_0])}.
$$
Thus,
\begin{equation*}
\begin{split}
| R_\epsilon(t) - R_0 | \leq 
C_3 T \Vert \tilde z_0 \Vert_{H_0^1([0,R_0])}, \quad \text{for } t \in [0,T]. 
\end{split}
\end{equation*}

If the take $K = C_2 \Vert \tilde z_0 \Vert_{L^\infty([0,R_0])}$ in \eqref{eq:espace_Z} and we assume  
$$\Vert \tilde z_0 \Vert_{H_0^1([0,R_0])} \leq min \left\lbrace \frac{K_1}{C_3}, \frac{R_0-R_*}{C_3 T}, \frac{E-R_0}{C_3 T}  \right\rbrace, $$
we have that
$\Vert R_\epsilon' \Vert_{L^\infty} \leq K_1$ and $R_* \leq R_\epsilon(t) \leq E$.
Thus $\Lambda_\epsilon(\mathcal{Y} \times \mathcal{Z}) \subset \mathcal{Y} \times \mathcal{Z}$.

The rest of the proof follows the same argument as the main Theorem of \cite{EFC_Ivaldo-17}. First, we have that  $\Lambda_\epsilon$ is compact, i.e. maps bounded sets of $L^\infty(\hat Q_0) \times C^1([0,T))$ into bounded sets of
$K^{0,\kappa}(\hat Q_0)\times C^{1,\kappa}([0,T])$.  

In particular, there exists $\kappa>0$ such that $\tilde z_\epsilon \in K^{0,\kappa}(\hat Q_0)$ and
\begin{equation}
\label{eq:z_eps_bounded}
|\tilde z_\epsilon|_{\hat Q_0}^{0,\kappa} \leq C,
\end{equation}
for some $C>0$ depending only on $R_*$, $E$, $T$, $\kappa$ and $\|\tilde z_0\|_{H_0^1([0,R_0])}$.
Moreover, from \eqref{eq:VR_bound}, 
there exists a constant $\tilde C>0$, only depending on the previous data, such that
\begin{equation}
\label{eq:R_eps_bounded}
\|R_\epsilon\|_{C^{1,\kappa}} \leq \tilde C.
\end{equation}

The continuity of the map $(z,R) \mapsto \Lambda_\epsilon(z,R)$ follows also as in \cite{EFC_Ivaldo-17}.

Applying Schauder's fixed point theorem to the map
$\Lambda_\epsilon : \mathcal{Y} \times \mathcal{Z} \to \mathcal{Y} \times \mathcal{Z}$, 
we get a fixed point $(\tilde z_\epsilon,R_\epsilon)$, for any $\epsilon>0$. The triple $(\tilde z_\epsilon, \tilde v_\epsilon,R_\epsilon)$ is a solution of \eqref{eq:1D-PDE-aprox1} and \eqref{eq:1D-PDE-aprox2}, for any $\epsilon>0$, and by \eqref{eq:z_eps_bounded}, \eqref{eq:v_bounded} and  \eqref{eq:R_eps_bounded}, each component of the triple is uniformly bounded. Thus, as $\epsilon \to 0$ we can extract subsequences such that
\begin{equation*}
\begin{cases}
R_\epsilon \to R \quad &\text{strongly in} \quad C^{1}([0,T]), \\
\tilde v_\epsilon \rightharpoonup \tilde v \quad &\text{weakly in} \quad L^2([0,b) \times (0,T)),
\end{cases}
\end{equation*}
where the strong converge is obtained by the compactness of the inclusion $C^{1+\kappa} \hookrightarrow C^1$.
We deduce from \eqref{eq:1D-PDE-aprox2}, as $\epsilon \to 0$, that $\tilde z$, the state associated to the control $\tilde v$ and the path $R(t)$, verifies  that $\tilde z(r,T)=0$, for $r \in [0,R(T)]$, and we conclude that the system \eqref{eq:1D-PDE-change} is locally null controllable. In the original variables this implies that the system \eqref{eq:PDE-intro}-\eqref{eq:stefan_condition-intro} is locally null controllable.
\end{proof}

\section{Additional remarks and open questions}
\label{sec:additional_remarks}

\subsection{Boundary controllability in two-dimensions}
In this section, we consider the case where the system is defined in an annular domain in two dimensions. That is, given a function $R(t) \in C^1([0,T])$ such that $0< A <R_*\leq R(t) \leq E$, $t \in (0,T)$, consider the annular domain $Q_{A,R} = \{ (x,t) : x \in B_{R(t)} - B_A, \ t \in (0,T) \}$ with exterior boundary
$\Sigma_R = \{ (x,t) : x \in \partial B_{R(t)}, \ t \in (0,T) \}$.

Consider the following problem in $Q_{A,R}$ where the interior boundary $\partial B_A$ is fixed, acted on by a boundary control $h$, and the exterior boundary $\Sigma_R$ follows a free-boundary condition as was considered in the previous sections. As above, let $f : \R \to \R$ be a globally Lipschitz continuous function such that $f(0)=0$ and consider the system
\begin{equation}
\label{eq:PDE-annular}
\begin{cases}
 y_t - \Delta y  + f(y) = 0, & \text{in } Q_{A,R}, \\
y(x,t) = 0, & \text{on } \Sigma_R, \\
y(x,t) = h(t), & \text{on } \partial B_{\bar A} \times (0,T), \\
 y(\cdot,0) = y_0, & \text{in } B_{R(0)} - B_{\bar A},
\end{cases}
\end{equation}
with Stefan condition
\begin{equation}
\label{eq:stefan_condition-annular}
R'(t) = -\nabla y(x,t) \cdot \hat n, \quad x \in \partial B_{R(t)}.
\end{equation}

Using the radial symmetry, we get an equation applied to the radial function $z(r,t)=y(x,t)$, where $r=|x|$. In this case the 2D Laplacian is given by $\Delta z = z_{rr} + \frac{1}{r}z_r$. The system \eqref{eq:PDE-annular} becomes
\begin{equation}
\label{eq:PDE-annular-radial}
\begin{cases}
z_t - z_{rr} - \frac{1}{r} z_r + f(z)=0 & \text{in } \hat Q_{A,R}, \\
z(A,t) = h(t), \quad z(R(t),t)=0  & \text{on } (0,T), \\
z(\cdot,0) = z_0 & \text{in } (A,R(0)),
\end{cases}
\end{equation}
where $\hat Q_{A,R} = \{ (r,t) : r \in [A,R(t)], \ t \in (0,T) \}$. Moreover $y(\cdot,t) \in L^2(B_{R(t)}-B_A)$ is equivalent to $z(\cdot,t) \in L^2([A,R(t)], r \ dr)$, as was observed in Section \ref{sec:reduction} for $N=3$.

Since the interior boundary $\partial B_A$ is fixed and $A > 0$, the factor $\frac{1}{r}$ in \eqref{eq:PDE-annular-radial} is bounded and we do not need to perform the second change of variables used in Section \ref{sec:reduction}.

We use a standard method to pass from boundary control to interior control: we work in an extended interval $[A-\delta,R(t)]$ and add an interior control $\bar v$ acting in an open subset $\bar \omega \subset (A-\delta,A)$. Then, the boundary control is given by $h(t) = \bar u(A,t)$ where $\bar u$ denotes the solution given by the controllability result of \cite{EFC_Ivaldo-17} with interior control $\bar v$ in the extended interval. As a consequence, we get the boundary controllability:

\begin{proposition}
Let $f$ be a globally Lipschitz continuous function such that $f(0)=0$, and $T>0$, $A,E>0$. Assume that $0 < A < R_* < R_0 < E$. Then, there exists $\varepsilon > 0$ such that, if $\Vert y_0 \Vert_{H_0^1(B_{R_0})} \leq \varepsilon$, there is a solution $(y,R,h)$ of \eqref{eq:PDE-annular}-\eqref{eq:stefan_condition-annular}, with $R(0)=R_0$, such that
\begin{align*}
R \in C^1([0,T]), \quad R_* \leq R(t) \leq E, \quad
h \in L^2([0,T]), \quad y \in C^0([0,T],H_0^1(B_{R(t)}),
\end{align*}
and $y(x,T) = 0$, for $x \in B_{R(T)}.$
\end{proposition}

We note that the $N$-dimensional generalization of this result holds, but in applications it is difficult to imagine a scenario where one acts on the interior boundary of an $N$-dimensional annular region for $N>2$.

\subsection{Related questions}

Next, we present some related open questions. We consider the same notation used throughout the paper.

\begin{itemize}
\item Consider the problem of the null controllability of the free-boundary  radial problem with general diffusion coefficient $C(|x|,t)$ given by
\begin{equation}
\label{eq:radial_PDE-open1}
\begin{cases}
y_t - div_x ( C(|x|,t) \ \nabla_x y ) + f(y) = 1_\omega u, & \text{in } Q_R, \\
y(x,t) = 0, & \text{on } \Sigma_R, \\
y(\cdot,0) = y_0, & \text{in } B_{R(0)},
\end{cases}
\end{equation}
with Stefan condition describing the evolution of the boundary,
\begin{equation*}
R'(t) = - C(|x|,t) \ \nabla y(x,t) \cdot \hat n, \quad x \in \partial B_{R(t)},
\end{equation*}
where $\hat n$ denotes the outward unit normal vector at the point $(x,t)$. We observe that we cannot use the same change of variable that we used in this paper, since depending on the form the function $C$, there is no certainty that we can cancel the singular term at the origin as in \eqref{eq:1D-PDE-change}.

\item Consider the function $\beta : \R \to \R$ of class $C^1$ with bounded derivatives and such that
$
0 < \beta_0 < \beta(r) < \beta_1 < +\infty.
$
Let us consider the problem of null controllability of the free-boundary radial problem with non-local term given by
\begin{equation}
\label{eq:radial_PDE-open2}
\begin{cases}
y_t - \beta\left( \int_0^{R(t)} y dx  \right) \Delta y + f(y) = 1_\omega u, & \text{in } Q_R, \\
y(x,t) = 0, & \text{on } \Sigma_R, \\
y(\cdot,0) = y_0, & \text{in } B_{R(0)},
\end{cases}
\end{equation}
with Stefan condition
\begin{equation*}
R'(t) = - \beta\left( \int_0^{R(t)} y dx  \right) \nabla y(x,t) \cdot \hat n, \quad x \in \partial B_{R(t)}.
\end{equation*}

The one-dimensional case has been studied in \cite{Vitor_JAL-23}  where higher regularity of the initial state was needed, i.e. $y_0 \in C^{2+\frac{1}{2}}([0,R_0])$, and also Hölder regularity of the approximate controls.

\item Let us assume that the real function $a = a(r)$ is of class $C^2(\R)$ and possesses bounded derivatives up to order two and such that
$
0 < C_0 < a(r) < C_1.
$
Consider the null controllability of the free-boundary quasi-linear radial problem given by
\begin{equation}
\label{eq:radial_PDE-open3}
\begin{cases}
y_t - div_x ( a(y) \nabla_x y ) + f(y) = 1_\omega u, & \text{in } Q_R, \\
y(x,t) = 0, & \text{on } \Sigma_R, \\
y(\cdot,0) = y_0, & \text{in } B_{R(0)},
\end{cases}
\end{equation}
with Stefan condition
\begin{equation*}
R'(t) = - a(y) \nabla y(x,t) \cdot \hat n, \quad x \in \partial B_{R(t)}.
\end{equation*}

The one-dimensional case was studied in \cite{WLL-22} based on the previous work \cite{Liu_Zhang-12}, where Hölder regularity was needed for the initial state and the controls obtained where also Hölder. In fact, $y_0 \in C^{2+\kappa}([0,R_0])$ and $u \in K^{\kappa,\frac{\kappa}{2}}$ for any $\kappa \in (0,1)$.

\end{itemize}

\backmatter





\bmhead{Acknowledgments}

This study was financed in part by the Coordenação de Aperfeiçoamento de Pessoal de Nível Superior - Brasil (CAPES) - Finance Code 001.
The authors thank Prof. Enrique Fernández-Cara for useful conversations on the Stefan problem and Martin Hernández and Daniel Fernández for useful comments on an earlier version of the manuscript. J. L. was partially supported by CNPq-Brazil. L. Y. was partially supported by CAPES-Brazil. 



\section*{Declarations}
\noindent Competing Interests: The authors have not disclosed any competing interests.









\begin{appendices}

\section{Proof of Proposition \ref{prop:carleman}}
\label{appendix A}

\begin{proof}
We follow \cite{EFC_Lim_Men-16} performing the computations on the interval $[0,R(t)]$ with control region $[0,b)$. 
Let us define $\psi$ by $\phi=e^{s\alpha} \psi$. Then
$$
\phi_t = s \alpha_t e^{s \alpha} \psi + e^{s \alpha} \psi_t, \qquad \phi_r = s \alpha_r e^{s \alpha} \psi + e^{s \alpha} \psi_r,$$
$$
\phi_{rr} = e^{s \alpha} \left( \psi_{rr} + s^2 \alpha_r^2 \psi + 2s\alpha_r \psi_r + s \alpha_{rr} \psi \right).
$$
Using the following formulas, for $k \geq 2$, 
$$
\alpha_r = -\xi_r = - \lambda \alpha_{0,r} \xi,
\qquad
\alpha_{rr} = -\lambda^2 \alpha_{0,r}^2 \xi - \lambda \alpha_{0,rr} \xi,
$$
$$
\alpha_t = - k \frac{T-2t}{t^{k+1}(T-t)^{k+1}} \sigma - \lambda \alpha_{0,t} \xi,
$$
and
$$
\alpha_{rt} = -\lambda^2 \xi \alpha_{0,t} \alpha_{0,x} - \lambda \xi \alpha_{0,rt} + k \lambda \xi \frac{T-2t}{t(T-t)} \alpha_{0,r},
$$
we have the estimates, for $\lambda$ sufficiently large, 
$$
|\alpha_r| \leq C \lambda \xi, \qquad |\alpha_{rr}| \leq C \lambda^2 \xi, \qquad |\alpha_t| \leq C \lambda \xi^{1+1/k} + C\lambda \xi \leq C\lambda\xi^{3/2}, \qquad |\alpha_{rt}|\leq C\lambda^2\xi.
$$
The adjoint equation \eqref{eq:1D-PDE-adjoint} becomes
$$
U\psi + V \psi = -e^{-\alpha}F - s\alpha_{rr}\psi - s\alpha_t \psi,
$$
where
$$
U\psi = \psi_t + 2s\alpha_r \psi_r, \qquad V\psi = \psi_{rr} + s^2 \alpha_r^2 \psi.
$$
Thus,
$$
\|U\psi\|_{L^2(\hat Q_R)}^2 + \|V\psi\|_{L^2(\hat Q_R)}^2 + 2(U\psi,V\psi)_{L^2(\hat Q_R)} = \|e^{-\alpha}F + s\alpha_{rr}\psi + s\alpha_t \psi\|_{L^2(\hat Q_R)}^2.
$$
We have
\begin{equation}
\label{eq:U-V}
(U\psi,V\psi)_{L^2(\hat Q_R)} = \int\int_{\hat Q_R} \left(\psi_t \psi_{rr} + s\alpha_r^3 \psi_t \psi + 2s\alpha_r\psi_r\psi_{rr} + 2s^3\alpha_r^3\psi_r\psi \right) drdt.
\end{equation}
The first term gives
\begin{equation*}
\begin{split}
\iint_{\hat Q_R} \psi_t \psi_{rr} &= -\frac{1}{2} \int_0^T\int_0^{R(t)} \frac{d}{dt}|\psi_r|^2 = -\frac{1}{2}\int_0^T \left( \frac{d}{dt}\int_0^{R(t)} |\psi_r|^2 dr - R'(t)|\psi_r(R(t),t)|^2 \right) dt \\
&=-\frac{1}{2}\int_0^{R(T)} |\psi_r(r,T)|^2 dr + \frac{1}{2} \int_0^{R(0)}|\psi_r(r,0)|^2 dr + \frac{1}{2}\int_0^T R'(t) |\psi_r(R(t),t)|^2 dt \\
&= \frac{1}{2}\int_0^T R'(t) |\psi_r(R(t),t)|^2 dt,
\end{split}
\end{equation*}
where we used that $\psi_r(r,0) \equiv 0$ in $(0,R(0))$ and $\psi_r(x,T) \equiv 0$ in $(0,R(T))$, by the definition of $\psi = e^{-s\alpha} \phi$ and the weight $\alpha$.

The third term of \eqref{eq:U-V} is computed by integration by parts paying attention to the boundary term, i.e.,
\begin{equation*}
\begin{split}
\iint_{\hat Q_R} 2s\alpha_r\psi_r\psi_{rr} &= -\iint_{\hat Q_R} s \alpha_{rr} |\psi_r|^2 + \int_0^T \left[ s \alpha_r |\psi_r|^2 \right]_{r=0}^{r=R(t)} \\
&= -\iint_{\hat Q_R} s \alpha_{rr} |\psi_r|^2 + \int_0^T \left( s \alpha_r(R(t),t) |\psi_r(R(t),t)|^2 - s \alpha_r(0,t) |\psi_r(0,t)|^2 \right) dt,
\end{split}
\end{equation*}
but $\alpha_r(0,t) = -\lambda \alpha_{0,r}(0,t) \xi(0,t)=0$, by the definition of the weight $\alpha_0$

The second and fourth terms of \eqref{eq:U-V} are computed by integration by parts where the boundary terms vanish by the boundary conditions of the adjoint equation \eqref{eq:1D-PDE-adjoint},
$$
\iint_{\hat Q_R} \left( s\alpha_r^3 \psi_t \psi + 2s^3\alpha_r^3\psi_r\psi \right) drdt = - \iint_{\hat Q_R} \left( s^2 \alpha_r \alpha_{rt} |\psi|^2 + 3s^3 \alpha_r^2 \alpha_{rr} |\psi|^2 \right) drdt.
$$

Thus \eqref{eq:U-V} becomes
\begin{equation*}
\begin{split}
(U\psi,V\psi)_{L^2(\hat Q_R)} =& s \iint_{\hat Q_R} (-\alpha_{rr}) |\psi_r|^2 drdt - \iint_{\hat Q_R} \left( 3s^3 \alpha_r^2 \alpha_{rr} + s^2 \alpha_r \alpha_{rt} \right)|\psi|^2 dr dt \\
&+ \int_0^T \left( s \alpha_r(R(t),t) + \frac{1}{2} R'(t) \right) |\psi_r(R(t),t)|^2 dt.
\end{split}
\end{equation*}
We use that
$$
\alpha_r(R(t),t) \geq C > 0 \quad \text{in} \quad (0,T),
$$
and
$$
-\alpha_r^2 \alpha_{rr} = \lambda^4 |\alpha_{0,r}|^4 \xi^3 + \lambda^3 |\alpha_{0,r}|^2\alpha_{0,rr} \xi^3 \geq C \lambda^4 \xi^3 - C \lambda^3 \xi^3, \qquad \text{in } \ \bar{\hat Q}_R - [0,b_0) \times (0,T),
$$
$$
|\alpha_r \alpha_{rt}| \leq C \lambda^3 \xi ^2, \qquad \text{in } \ \bar{\hat Q}_R - [0,b_0) \times (0,T).
$$
Taking $\lambda$ sufficiently large, 
\begin{equation*}
\begin{split}
(U\psi,V\psi)_{L^2(\hat Q_R)} \geq & s \iint_{\hat Q_R} (-\alpha_{rr}) |\psi_r|^2 drdt + \iint_{\hat Q_R - [0,b_0) \times (0,T)} \lambda^4 s^3 \xi^3 dr dt \\
& + \ s \lambda \int_0^T \xi(R(t),t) |\psi_r(R(t),t)|^2 dt.
\end{split}
\end{equation*}

Thus, from the estimate
$$\|e^{-\alpha}F + s\alpha_{rr}\psi + s\alpha_t \psi\|_{L^2(\hat Q_R)}^2 \leq C \iint_{\hat Q_R} e^{-2\alpha}|F|^2drdt + C \iint_{\hat Q_R} \lambda^4 s^2 \xi^3 |\psi|^2 drdt,$$
we have
\begin{equation}
\label{eq:A2}
\begin{split}
\|U\psi\|_{L^2(\hat Q_R)}^2 + \|V\psi\|_{L^2(\hat Q_R)}^2 + s\lambda\int_0^T \xi(R(t),t)|\psi_r(R(t),t)|^2 dt \\
+ \iint_{\hat Q_R - [0,b_0) \times (0,T)} \left( \lambda^2 (s\xi)|\psi_r|^2 + \lambda^4 (s\xi)^3 |\psi|^2 \right)drdt \\ 
\leq C \iint_{\hat Q_R} e^{-2\alpha}|F|^2drdt + C \iint_{\hat Q_R} \lambda^4 s^2 \xi^3 |\psi|^2 drdt. 
\end{split}
\end{equation}

Taking $s$ sufficiently large to absorb the term with $\lambda^4 s^2$ in \eqref{eq:A2}, 
\begin{equation}\label{eq:A3}
\begin{split}
&\|U\psi\|_{L^2(\hat Q_R)}^2 + \|V\psi\|_{L^2(\hat Q_R)}^2 + s\lambda\int_0^T \xi(R(t),t)|\psi_r(R(t),t)|^2 dt \\
&+ \iint_{\hat Q_R} \left( \lambda^2 (s\xi)|\psi_r|^2 + \lambda^4 (s\xi)^3 |\psi|^2 \right)drdt \\
&\leq C \iint_{\hat Q_R} e^{-2\alpha}|F|^2drdt + C \iint_{[0,b_0) \times (0,T)} \lambda^4 (s \xi)^3 |\psi|^2 drdt +
C \iint_{[0,b_0) \times (0,T)} \lambda^2 (s \xi) |\psi_r|^2 drdt.
\end{split}
\end{equation}

Arguing as in \cite{Fur_Ima-96}, we can absorb the integral with $|\psi_r|^2$ in the right-hand side of \eqref{eq:A3}. On the other hand, we have
\begin{gather*}
\|U \psi\|_{L^2(\hat Q_R)}^2 \geq C\iint_{\hat Q_R} (s\xi)^{-1}|\psi_{t}|^2 - C \iint_{\hat Q_R} \lambda^2 (s\xi) |\psi_r|^2, \\
\|V \psi\|_{L^2(\hat Q_R)}^2 \geq C\iint_{\hat Q_R} (s\xi)^{-1}|\psi_{rr}|^2 - C \iint_{\hat Q_R} \lambda^4 (s\xi)^3 |\psi|^2.
\end{gather*}

Finally, coming back to the original variable $\phi$ we conclude \eqref{eq:carleman}.
\end{proof}




\end{appendices}



\end{document}